\documentclass[11pt]{amsart}
\usepackage{amssymb}
\usepackage{fullpage}
\usepackage{xcolor}
\usepackage[hidelinks]{hyperref}
\usepackage[capitalise]{cleveref}

\numberwithin{equation}{section}

\newtheorem*{corollary*}{Corollary}
\newtheorem*{theorem*}{Theorem}
\newtheorem{theorem}{Theorem}[section]

\newtheorem{lemma}[theorem]{Lemma}

\theoremstyle{definition}

\newcommand{\E}{\mathbb E}
\newcommand{\C}{\mathbb C}
\DeclareMathOperator{\tr}{tr}
\DeclareMathOperator{\Sym}{Sym}

\title{Harper's beyond square-root conjecture for random matrices}
\date{}

\author{William Verreault}
\address{Department of Mathematics \\
University of Toronto \\
Toronto, ON \\
Canada}
\email{william.verreault@utoronto.ca}

\author{Max Wenqiang Xu}
\address{Yau Mathematical Sciences Center, Tsinghua University, Beijing, 100084, China}
\address{Beijing Institute of Mathematical Sciences and Applications, Beijing, 101408, China}
\email{maxxu1729@gmail.com, maxxu@tsinghua.edu.cn}

\begin{document}

\vspace*{-4mm}

\begin{abstract}
For a Haar unitary matrix $U_N$ of size $N$, we prove that the first absolute moment of $\tr(\Sym^k U_N)$ tends to zero as $\min(N,k)\to\infty$. This settles the random matrix version of Harper's beyond square-root conjecture formulated by Wang and Xu in a stronger form, as our result does not have any restriction on $k/N$. The proof is elementary and self-contained, combining exact weighted trace moments with a Gaussian change of measure. It connects the result to the critical chaos model studied by Soundararajan and Zaman, also known as holomorphic multiplicative chaos.
\end{abstract}

\maketitle

\section{Introduction} 
\subsection{Background}
Harper's theorem on ``better than square-root cancellation'' in partial sums of random multiplicative functions \cite{HarperRMF} marked a major advance in probabilistic number theory, revealing a surprising connection with critical multiplicative chaos. We refer to \cite{HarperReview,HarperSurvey} for surveys of these ideas and their number theoretic significance.

A particularly attractive feature of Harper's theorem is that it also sheds light on fundamental deterministic problems in number theory.
For example, using  his derandomization argument, Harper \cite{HarperCharacter} proved that for a prime modulus $r$,
\[ \frac{1}{r-1} \sum_{\chi\bmod r} |\sum_{n\le x} \chi(n)| = o(\sqrt{x})\]
whenever $\min(x,r/x)\to\infty$. An analogous result holds for the continuous family $n^{it}$, averaged over $|t|\le T$, when $\min(x,T/x)\to\infty$. A further application to Möbius cancellation is particularly striking. Harper proved that the same phenomenon holds in the same range after twisting by the Möbius function $\mu$,
and he conjectured that this should hold for $x\le r^{100}$. It turns out that this seemingly unimportant summation range $x\le r^{100}$ is crucial. Indeed, it is not hard to show that if Harper's conjecture holds, then Möbius cancellation holds in the form
\[\sum_{x\le n \le x+y} \mu(n) = o(y) \]
as long as $y \gg \sqrt{x}/W(x)$ for some $W(x)$ tending to infinity. In fact, Harper's conjecture is quantitative, and as a corollary, one can take $W(x) = (\log \log x)^{c}$ for any fixed $0<c<1/4$. This would break the famous square-root barrier,
and note that assuming the Riemann Hypothesis can only push the range of $y$ to $\gg x^{1/2+\epsilon}$. See also \cite{HarperCharacter,HarperSurvey}. 
In \cite{WX}, Wang and Xu managed to prove a version of Harper's conjecture under strong assumptions, namely the Ratios conjecture and the Generalized Riemann Hypothesis.

In this paper, we establish Harper's conjecture unconditionally in the random matrix setting, which one may also view as a function field analogue. Since Möbius twists occur as coefficients in the Dirichlet series of reciprocals of $L$-functions, the relevant random matrix quantities are the coefficients
\begin{equation}\label{eq:coeff}
h_k(U)=[z^k]\det(I-zU)^{-1}=\tr(\Sym^k U),
\qquad U\in U(N).
\end{equation}
Wang and Xu conjectured that the corresponding first absolute moment tends to $0$ as $k\to\infty$, provided $k\leq BN$ for a fixed $B>0$ \cite[Conjecture 1.5]{WX}. In the arithmetic setting, $k$ corresponds to $\log x$ and $N$ to the logarithm of the conductor, so a bounded ratio $k/N$ corresponds to sums of length at most a fixed power of the conductor. We prove their conjecture in a stronger form, i.e., with no restriction on $k/N$.

\subsection{Main result}
We write $\E_{U(N)}$ for integration against the normalized Haar measure on $U(N)$.
\begin{theorem}\label{thm:main}
Let $U_N$ be Haar distributed on $U(N)$. Then
$$
\E_{U(N)}\big|\tr(\Sym^k U_N)\big|\longrightarrow0 \qquad\text{as }\min(N,k)\longrightarrow\infty.
$$
\end{theorem}
We remark that the method here can likely be adapted to other compact groups, e.g., the symplectic and orthogonal families and we refer readers to \cite{HamdanLeungWong2026LowMoments, XuZheng2026TypicalSize} for related results on other families. 

\subsection{Connections with function fields and multiplicative chaos}
We give a quick explanation to show the relation with function fields (see also \cite{SZ}). For a primitive odd Dirichlet character $\chi$ modulo a monic squarefree polynomial $Q\in\mathbb F_q[T]$ (with $\deg Q\ge2$), write $L(u,\chi)=\det(I-\sqrt q\,u\,\Theta_\chi)$, where $\Theta_\chi \in U(\deg Q-1)$. Then one has 
$$
q^{-k/2}\sum_{\substack{f\ \mathrm{monic}\\ \deg f=k}}\mu(f)\chi(f)
=h_k(\Theta_\chi).
$$
Katz's equidistribution theorem \cite{KatzSquarefree} shows that averages over these characters converge to Haar averages as $q\to\infty$, with $k$ and $\deg Q$ fixed. This gives a function field interpretation of our theorem, with the field size tending to infinity first. See also related work of Gorodetsky \cite{Gorodetsky} and the papers \cite{fusrp2,fusrp1} on moment estimates and distribution theorems for random multiplicative functions over function fields. For geometric estimates for sums of arithmetic functions in short intervals and squarefree arithmetic progressions over function fields, see Sawin \cite{SawinDuke,SawinActa}.

The relation between random multiplicative functions and critical chaos has been instrumental in many recent works. See \cite[Section 1.6]{Verreault} for a recent account. 
This becomes especially transparent through the multiplicative chaos model studied by Soundararajan and Zaman \cite{SZ}. Let $X(j)$ be independent standard complex Gaussian random variables. We write $\E$ for expectation over these and any other auxiliary variables introduced below. Define $A(k)$ by
\begin{equation}\label{eq:gaussian-series}
\sum_{k\ge0}A(k)z^k
=\exp\Big(\sum_{j\ge1}\frac{X(j)}{\sqrt j}z^j\Big).
\end{equation}
The random power series in \eqref{eq:gaussian-series} is also studied as critical holomorphic multiplicative chaos by Najnudel, Paquette and Simm \cite{NPS}.
These coefficients arise as large-field limits of normalized sums of Steinhaus random multiplicative functions over polynomials of fixed degree \cite{SZ}. For fixed $k$, they also arise as the limits in distribution of $h_k(U_N)$ as $N\to\infty$, by the trace limit theorem of Diaconis and Shahshahani \cite{DiaconisShahshahani}. Soundararajan and Zaman proved that
$$
\E|A(k)|^2=1,
\qquad
\E|A(k)|\asymp(\log k)^{-1/4},
$$
which can be viewed as a function field version of Harper's theorem. An asymptotic refinement of this estimate was conjectured and supported numerically in \cite{ASVZZ}. 
Together with the second moment identity, the fixed $k$ convergence gives the iterated limit $$\lim_{k\to\infty}\lim_{N\to\infty}\E_{U(N)}|h_k(U_N)|=0.$$ It does not, however, provide the uniformity needed when $k$ grows with $N$. 

For the ordinary characteristic polynomial, whose coefficient of $z^k$ is $(-1)^k\tr(\wedge^k U_N)$, with $\wedge^k$ denoting the $k$th exterior power, Najnudel, Paquette and Simm \cite[Theorem~1.2]{NPS} proved convergence to zero in distribution as $k\to\infty$ with $k\le N/2$. Since the second moment equals one, the first absolute moment also tends to zero. In comparison, our theorem concerns the symmetric powers arising from the reciprocal characteristic polynomial and allows arbitrary relative growth of $k$ and $N$. We remark that this type of coefficient problem is also related to the harmonic analysis of Gaussian multiplicative chaos measures (see \cite{OrsoniVerreault}).

\section{Notation and preliminaries}

\subsection{Partitions and symmetric functions}
For a nonincreasing sequence of positive integers $\lambda_i$, we use $\lambda=(\lambda_1,\ldots,\lambda_\ell)$ to denote a partition of length $\ell(\lambda)=\ell$, and $|\lambda|=\sum_i\lambda_i$ for its size. Define $m_j(\lambda)=\#\{i:\lambda_i=j\}$, the number of parts in $\lambda$ equal to $j$. 
We write $\lambda\vdash n$ if $|\lambda|=n$. Put
$$
z_\lambda=\prod_{j\ge1}j^{m_j(\lambda)}m_j(\lambda)!.
$$
We define the power sums $p_j$ and the complete homogeneous symmetric polynomials $h_n$ as $$p_j(x)=\sum_{i=1}^N x_i^j \quad \text{and} \quad
h_n(x)=\sum_{1\le i_1\le\cdots\le i_n\le N}x_{i_1}\cdots x_{i_n}.
$$
Recall the classical generating function identity 
\begin{equation}\label{eq:complete}
\sum_{n\ge0}h_n(x)z^n =\prod_{i=1}^N(1-x_i z)^{-1} =\exp\Big(\sum_{j\ge1}\frac{p_j(x)z^j}{j}\Big).
\end{equation}
Expanding the exponential and grouping terms by the multiplicities gives
$$
h_n=\sum_{\lambda\vdash n}\frac{p_\lambda}{z_\lambda}, \qquad p_\lambda=\prod_{j\geq 1}p_j^{m_j(\lambda)}.
$$
We state its connection to random matrix statistics. If the $x_i$ are the eigenvalues of $U\in U(N)$, then $p_j(U)=\tr U^j$ and $h_n(U)=\tr(\Sym^n U)$ as in \eqref{eq:coeff}.

\subsection{Inner product structure}
Let $X(j)$ be independent standard complex Gaussian random variables. We use the notation $\mathbf g=(\sqrt j X(j))_{j\ge1}$, then we have
$$
h_n(\mathbf g)=A(n)=\sum_{\lambda\vdash n} \prod_{j\ge1}\frac{X(j)^{m_j(\lambda)}}{j^{m_j(\lambda)/2}m_j(\lambda)!}.
$$
Since $\E[X(j)^a\overline{X(j)}^{b}]=\mathbf{1}_{a=b}a!$, independence gives
$
\E[p_\lambda(\mathbf g)\overline{p_\rho(\mathbf g)}] = \mathbf{1}_{\lambda=\rho}\,z_\lambda.
$ Writing $P(\mathbf g)$ for the specialization $p_j=\sqrt j X(j)$, we have an inner product
\begin{equation} \label{eq:Hall}
\langle P,Q\rangle=\E[P(\mathbf g)\overline{Q(\mathbf g)}].
\end{equation}

\subsection{Schur polynomials}
We recall some classical facts about symmetric polynomials and we refer readers to \cite[Chapter I]{Mac}.

For a partition $\lambda$ of length $\ell$, define its Schur function by the determinant
$$
s_\lambda=\det\big(h_{\lambda_i-i+j}\big)_{1\le i,j\le\ell}.
$$
In particular $s_{(k)}=h_k$ for $k\geq 1$. 
For each $n$, the $s_\lambda$ with $\lambda\vdash n$ form a basis of the homogeneous polynomials of weighted degree $n$ in the $p_j$ and are orthonormal for the inner product \eqref{eq:Hall}. In $N$ variables, $s_\lambda$ vanishes if $\ell(\lambda)>N$. For $\ell(\lambda)\le N$, $s_\lambda(U)$ is an irreducible character of $U(N)$, so the specializations of these Schur functions are also orthonormal under Haar probability.

The Young diagram of $\lambda$ has $\lambda_i$ boxes in row $i$. Pieri's rule states that
\begin{equation}\label{eq:pieri}
h_k s_\lambda=\sum_{\rho:\,\rho/\lambda\text{ is a horizontal }k\text{-strip}}s_\rho,
\end{equation}
where a horizontal $k$-strip means that $\rho$ is obtained from $\lambda$ by adding $k$ boxes, with no two added boxes in the same column. Equivalently,
$$
\rho_1\ge\lambda_1\ge\rho_2\ge\lambda_2\ge\cdots, \qquad |\rho|-|\lambda|=k.
$$
A consequence of this is that $\ell(\rho)\le\ell(\lambda)+1.$  This bound is independent of the number $k$ of boxes added.

\section{Proof of \cref{thm:main}} \label{sec:proof}

\subsection{Exact moment comparison}
We will use a lemma that allows us to compare weighted moments. The identity below follows from a moment formula for complete symmetric functions, for instance as in recent work \cite[Section 3]{Gorodetsky}. It can also be deduced by coefficient extraction from the unitary ratios formula of Conrey, Farmer, and Zirnbauer \cite{CFZ2005}. We give a short proof for completeness.

\begin{lemma}\label{lem:transfer}
Let $P,Q$ be polynomials in the power sums of weighted degree at most $D$. For every $k\ge0$ and $N\ge D+1$,
$$
\E_{U(N)}\big[|h_k(U)|^2P(p(U))\overline{Q(p(U))}\big] =\E\big[|A(k)|^2P(\mathbf g) \overline{Q(\mathbf g)}\big].
$$ 
Here, $p(U)$ stands for $(p_j(U))_{j\ge1}$, $A(k)$ is defined in \eqref{eq:gaussian-series} and $\mathbf g=(\sqrt j X(j))_{j\ge1}$.
\end{lemma}

\begin{proof}
Expand $P$ and $Q$ in the Schur basis. Note that only partitions of size at most $D$ occur, and these have length at most $D$.
By Pieri's rule \eqref{eq:pieri}, every Schur function occurring in $h_kP$ or $h_kQ$ therefore has length at most $D+1\le N$. These Schur functions are orthonormal for the inner product \eqref{eq:Hall}, and after evaluation at the eigenvalues of $U$, for the Haar inner product too. This leads to
$$
\E_{U(N)}\big[(h_kP)(U)\overline{(h_kQ)(U)}\big]
=\E\big[(h_kP)(\mathbf g)\overline{(h_kQ)(\mathbf g)}\big].
$$
Since $h_k(\mathbf g)=A(k)$, this completes the proof.
\end{proof}

The special case $P=Q=1$ gives $\E|A(k)|^2=\E_{U(N)}|h_k|^2=1$. At $k=0$, the comparison reduces to the classical trace moment identity \cite{DiaconisShahshahani}.

We next introduce suitable weighted measures.
Let $\mu_N$ be Haar probability and let $\nu_{N,k}$ be the probability measure defined by $d\nu_{N,k}=|h_k|^2d\mu_N$. Fix $m$ and write
\begin{equation}\label{eqn: notation}
T_m(U)=\big(p_j(U)/\sqrt j\big)_{1\le j\le m},\qquad X_m=(X(j))_{j=1}^m, \qquad W_m=\big(X(j)+e^{ij\Theta}/\sqrt j\big)_{1\le j\le m}, 
\end{equation}
where $\Theta$ is uniformly distributed on $[0,2\pi)$ and independent of the Gaussian sequence. Recall that $\E$ denotes expectation over all the auxiliary variables involved.

\begin{lemma}\label{lem:weighted}
Let $P(z,\overline z)$ be a polynomial in $z_1,\ldots,z_m$ and their conjugates, such that every monomial $z^\alpha\overline z^\beta$ occurring in $P$ satisfies degree conditions
$$
\sum_{j=1}^m j\alpha_j\le D,\qquad \sum_{j=1}^m j\beta_j\le D.
$$
Then for $N\geq D+1$ and $k\geq D$, we have
$$
\E_{\nu_{N,k}}P(T_m,\overline{T_m}) =\E P(W_m,\overline{W_m}).
$$
\end{lemma}

\begin{proof}
Using \eqref{eq:gaussian-series} and the complex Gaussian identity $\E e^{aX(j)+b\overline{X(j)}}=e^{ab}$, we get that
\begin{align*}
\sum_{a,b\ge0}x^ay^b\, \E\Big[A(a)\overline{A(b)} \exp\Big(\sum_{j\le m}(t_jX(j)+s_j\overline{X(j)})\Big)\Big] = \frac{\exp(\sum_{j\le m}t_js_j)}{1-xy} \exp\Big(\sum_{j\le m}\frac{t_jy^j+s_jx^j}{\sqrt j}\Big).
\end{align*}
Fix a monomial $t^\alpha s^\beta$ with $\sum_jj\alpha_j,\sum_jj\beta_j\leq D$. For every term $x^ay^b$ from the last exponential that contributes to this coefficient, we must have $a,b\leq D$. For $k\geq D$, we therefore have
$$
[x^ky^k]\frac{x^ay^b}{1-xy}=\mathbf{1}_{a=b}.
$$
For the chosen coefficient, extracting the coefficient of $x^ky^k$ consequently has the same effect as retaining the factor $\exp(\sum_{j\leq m}t_js_j)$ and averaging the last exponential with $x=e^{-i\theta}$ and $y=e^{i\theta}$. By using Fourier orthogonality, the resulting coefficient is that of $t^\alpha s^\beta$ in
$$
\exp\Big(\sum_{j\le m}t_js_j\Big) \int_0^{2\pi}\exp\Big(\sum_{j\le m} \frac{t_je^{ij\theta}+s_je^{-ij\theta}}{\sqrt j}\Big)\frac{d\theta}{2\pi}.
$$
Conditioning on $\Theta$ and applying the Gaussian identity again shows that the above is $$\E\exp\Big(\sum_{j\le m}(t_jW_{m,j}+s_j\overline{W_{m,j}})\Big).$$ 
Here $W_{m, j}$ is clearly the $j$-th coordinate of $W_m$. 

Comparing coefficients and multiplying by $\alpha!\beta!$ gives
$
\E\big[|A(k)|^2X_m^\alpha\overline{X_m}^{\beta}\big]
=\E\big[W_m^\alpha\overline{W_m}^{\beta}\big].
$
Applying \cref{lem:transfer} with the polynomials $\prod_{j\le m}(p_j/\sqrt j)^{\alpha_j}$ and $\prod_{j\le m}(p_j/\sqrt j)^{\beta_j}$ proves the identity for monomials. Hence the conclusion holds for $P$ by linearity.
\end{proof}

\subsection{Completion of the proof of \cref{thm:main}}
Let $\gamma_m$ be the law of Gaussian vector $X_m$, write $H_m=\sum_{j\le m}j^{-1}$ for the harmonic sum, and define
\begin{equation}\label{eq:mass}
M_m(z)=e^{-H_m}\int_0^{2\pi} \exp\Big(2\Re\sum_{j\le m}\frac{z_je^{ij\theta}}{\sqrt j}\Big) \frac{d\theta}{2\pi}.
\end{equation}
The density of a fixed translate $X_m+v$ relative to $\gamma_m$ is $\exp(2\Re\sum_jz_j\overline v_j-\sum_j|v_j|^2)$. 
Averaging this over $v_j=e^{ij\Theta}/\sqrt j$ and changing $\theta$ to $-\theta$ in the integral shows that $W_m$ has density $M_m$ relative to $\gamma_m$. It is clear that $M_m$ is positive and continuous, and $\E M_m(X_m)=1$. 

For fixed $m$, \cref{lem:transfer} with $k=0$ shows that every fixed mixed moment of $T_m$ (as defined in \eqref{eqn: notation}) under $\mu_N$ agrees with the corresponding moment of $X_m$ for all sufficiently large $N$. By \cref{lem:weighted}, every fixed mixed moment under $\nu_{N,k}$ agrees with the corresponding moment of $W_m$ for all sufficiently large $N$ and $k$. The second moments imply tightness, while higher even moments imply uniform integrability of every fixed polynomial. Consequently, every subsequential limit has the moments of $X_m$ or $W_m$, respectively. It is also clear that these laws are determined by their moments.

Thus, for every bounded continuous $f:\C^m\to\C$, as $\min(N,k)\to \infty$,
\begin{equation}\label{eq:limits}
    \begin{aligned}
\E_{U(N)}f(T_m) \longrightarrow \E f(X_m), \quad \E_{\nu_{N,k}}f(T_m) \longrightarrow \E[M_m(X_m)f(X_m)].
    \end{aligned}
\end{equation}

We next establish the qualitative estimate needed to complete the proof. A quantitative estimate can also be established, and we remark on this in the next subsection.

\begin{lemma} \label{lem:critical-mass}
Let the notation be the same as before. Then as $m\to\infty$, 
$$
\E \frac{M_m(X_m)}{1+M_m(X_m)}\longrightarrow 0.
$$    
\end{lemma}

\begin{proof}
Set
\[R_m=M_m(W_m)^{-1}.\] Since $W_m$ has density $M_m$ relative to $\gamma_m$, for $n\geq m$ and every bounded continuous $f:\C^m\to \C$, we have
$$
\E[f(W_m)R_n]
= \E f(X_m)
= \E[f(W_m)R_m].
$$
Thus $(R_m)$ is a nonnegative martingale of mean one for the filtration $(\sigma(W_m))_{m\geq 1}$, and consequently it converges almost surely to a finite limit.

We claim that $\limsup_m M_m(W_m)=\infty$ almost surely. By rotation invariance, the sequence $(M_m(W_m))_{m\geq 1}$ has the same law as 
$$
\big(M_m\big((X(j)+j^{-1/2})_{j\leq m}\big)\big)_{m\geq 1}.
$$
Set
$$
a_{j,m}=\frac{\sin(j/m)}{j/m}, \qquad G_m=2\Re\sum_{j\leq m}\frac{a_{j,m}X(j)}{\sqrt{j}}.
$$
Using periodicity, restricting the integral in \eqref{eq:mass} to $[-1/m,1/m]$ and applying Jensen's inequality gives 
$$
M_m\big((X(j)+j^{-1/2})_{j\leq m}\big)
\geq \frac{1}{\pi m}e^{-H_m+2\sum_{j\leq m}\frac{a_{j,m}}{j}+G_m} \geq ce^{G_m}
$$
for some absolute constant $c>0$. Here we used $a_{j,m}=1+O((j/m)^2)$ and $H_m=\log m+O(1)$.
The variable $G_m$ is a centered Gaussian with
$$
\operatorname{Var}(G_m)=2\sum_{j\leq m}\frac{a_{j,m}^2}{j} = 2\log m + O(1).
$$
Hence, for each real $A$, $\mathbb P(G_m>A)\to 1/2$, which implies $\mathbb P(\sup_m G_m=\infty) \geq 1/2$. It is clear that changing finitely many of the variables $X(j)$ changes the sequence $(G_m)$ by an amount bounded uniformly in $m$. Thus $\{\sup_m G_m=\infty\}$ is a tail event, and Kolmogorov's zero-one law gives $\mathbb P(\sup_m G_m=\infty)=1$. This proves the claim.
It follows that $\liminf_m R_m=0$ almost surely. Since $(R_m)$ converges, we have $R_m\to 0$ almost surely. Finally, using the definition of $R_m$ and the bounded convergence give
$$
\E\frac{M_m(X_m)}{1+M_m(X_m)} = \E\frac{1}{1+M_m(W_m)} = \E\frac{R_m}{1+R_m} \longrightarrow 0,
$$
as desired.
\end{proof}

Set $g_m=(1+M_m)^{-1}$ and $g=g_m(T_m)$. The Cauchy--Schwarz inequality gives
$$
\E_{U(N)}|h_k| =\E_{U(N)}(|h_k|g)+\E_{U(N)}(|h_k|(1-g)) \le\sqrt{\E_{\nu_{N,k}}g}+\sqrt{\E_{U(N)}(1-g)}.
$$
Here we used that $\mu_N$ and $\nu_{N,k}$ are probability measures, and $0\leq g\leq 1$.
By \eqref{eq:limits}, both expectations on the right tend to
$
\E [M_m(X_m)/(1+M_m(X_m))]$. Therefore
$$
\limsup_{\min(N,k)\to\infty}\E_{U(N)}|h_k(U_N)| \le2\Big(\E\frac{M_m(X_m)}{1+M_m(X_m)}\Big)^{1/2}.
$$
An application of \cref{lem:critical-mass} proves \cref{thm:main}.

\subsection{Connection with critical chaos} \label{sec:critical}
    The random variable $M_m(X_m)$ is the mean one total mass of the Fourier approximation to critical Gaussian multiplicative chaos on the circle. This gives the connection with the critical chaos model discussed in the introduction. 

    Given this, one can also deduce \cref{lem:critical-mass} directly from the low moment estimates for the total mass of the chaos, established by Soundararajan and Zaman. Indeed, as in \cite{SZ}, write
$
F_m(z)=\exp\big(\sum_{j\le m}\frac{X(j)z^j}{\sqrt j}\big).
$
By \eqref{eq:mass},
$$
M_m(X_m)=\frac{e^{-H_m}}{2\pi} \int_0^{2\pi}|F_m(e^{i\theta})|^2\,d\theta.
$$
The low moment estimates of \cite[Proposition 3.2]{SZ}, with $K=m$, $r=1$ and $q=1/2$, give $$\E\sqrt{M_m(X_m)}\ll (\log m)^{-1/4}\to0,$$ since $H_m=\log m+O(1)$. The inequality $x/(1+x)\leq \sqrt{x}$ then implies the desired result.

\subsection*{AI disclosure} AI tools were used to help with a literature search for moment formulas and to find the exact result of Conrey, Farmer, and Zirnbauer referred to in \cref{sec:proof}. They also helped with formatting the bib file for references, and with language review and polishing at the final stage of writing.

\subsection*{Acknowledgments} This work was initiated during the first author's visit to Tsinghua University, whose hospitality he gratefully acknowledges. We thank Victor Wang for helpful discussions and for pointing out the connection between \cref{lem:transfer} and the theorem of Conrey, Farmer, and Zirnbauer, as well as Ofir Gorodetsky for his helpful remarks on the function field analogue and the work of Sawin.

\bibliographystyle{amsplain}
\bibliography{references}

@article{DiaconisShahshahani,
  author       = {Diaconis, Persi and Shahshahani, Mehrdad},
  title        = {On the eigenvalues of random matrices},
  journal      = {J. Appl. Probab.},
  fjournal     = {Journal of Applied Probability},
  volume       = {31A},
  year         = {1994},
  pages        = {49--62},
  note         = {Studies in applied probability},
}

@article{Gorodetsky,
  author       = {Gorodetsky, Ofir},
  title        = {Magic squares, the symmetric group and {M}\"{o}bius randomness},
  journal      = {Monatsh. Math.},
  fjournal     = {Monatshefte f\"{u}r Mathematik},
  volume       = {204},
  number       = {1},
  year         = {2024},
  pages        = {27--46},
}

@article{HarperRMF,
  author       = {Harper, Adam J.},
  title        = {Moments of random multiplicative functions, {I}: {L}ow moments,
                  better than squareroot cancellation, and critical
                  multiplicative chaos},
  journal      = {Forum Math. Pi},
  fjournal     = {Forum of Mathematics. Pi},
  volume       = {8},
  year         = {2020},
  pages        = {e1, 95},
}

@incollection{HarperSurvey,
  author       = {Harper, Adam J.},
  title        = {Better than squareroot cancellation in number theory},
  booktitle    = {Proceedings of the {International Congress of Mathematicians}},
  publisher    = {SIAM},
  year         = {2026},
  pages        = {396--413},
}

@article{KatzSquarefree,
  author       = {Katz, Nicholas M.},
  title        = {On a question of {K}eating and {R}udnick about primitive
                  {D}irichlet characters with squarefree conductor},
  journal      = {Int. Math. Res. Not. IMRN},
  fjournal     = {International Mathematics Research Notices. IMRN},
  number       = {14},
  year         = {2013},
  pages        = {3221--3249},
}

@book{Mac,
  author       = {Macdonald, I. G.},
  title        = {Symmetric functions and {H}all polynomials},
  series       = {Oxford Mathematical Monographs},
  edition      = {Second},
  publisher    = {The Clarendon Press, Oxford University Press, New York},
  year         = {1995},
  pages        = {x+475},
  isbn         = {0-19-853489-2},
}

@article{NPS,
  author       = {Najnudel, Joseph and Paquette, Elliot and Simm, Nick},
  title        = {Secular coefficients and the holomorphic multiplicative chaos},
  journal      = {Ann. Probab.},
  fjournal     = {The Annals of Probability},
  volume       = {51},
  number       = {4},
  year         = {2023},
  pages        = {1193--1248},
}

@article{SawinDuke,
  author       = {Sawin, Will},
  title        = {Square-root cancellation for sums of factorization functions
                  over short intervals in function fields},
  journal      = {Duke Math. J.},
  fjournal     = {Duke Mathematical Journal},
  volume       = {170},
  number       = {5},
  year         = {2021},
  pages        = {997--1026},
}

@article{SZ,
  author       = {Soundararajan, Kannan and Zaman, Asif},
  title        = {A model problem for multiplicative chaos in number theory},
  journal      = {Enseign. Math.},
  fjournal     = {L'Enseignement Math\'{e}matique},
  volume       = {68},
  number       = {3-4},
  year         = {2022},
  pages        = {307--340},
}

@article{WX,
  author       = {Wang, Victor Y. and Xu, Max Wenqiang},
  title        = {{H}arper's beyond square-root conjecture},
  journal      = {Int. Math. Res. Not. IMRN},
  fjournal     = {International Mathematics Research Notices. IMRN},
  number       = {18},
  year         = {2025},
  pages        = {Paper No. rnaf279, 22},
}

@misc{CFZ2005,
  author       = {Conrey, J. B. and Farmer, D. W. and Zirnbauer, M. R.},
  title        = {{Howe} pairs, supersymmetry, and ratios of random
                  characteristic polynomials for the unitary groups {$U(N)$}},
  howpublished = {arXiv:math-ph/0511024},
  year         = {2005},
  url          = {https://arxiv.org/abs/math-ph/0511024},
}

@misc{HamdanLeungWong2026LowMoments,
  author       = {Hamdan, Jad and Leung, Sun-Kai and Wong, Mo Dick},
  title        = {Low moments of {Hecke} eigenvalue sums},
  howpublished = {arXiv:2609.15937},
  year         = {2026},
  url          = {https://arxiv.org/abs/2609.15937},
}

@misc{HarperCharacter,
  author       = {Harper, Adam J.},
  title        = {The typical size of character and zeta sums is {$o(\sqrt{x})$}},
  howpublished = {arXiv:2301.04390},
  year         = {2023},
  url          = {https://arxiv.org/abs/2301.04390},
}

@misc{HarperReview,
  author       = {Harper, Adam J.},
  title        = {Moments of random multiplicative functions, {III}: {A} short
                  review},
  howpublished = {arXiv:2410.11523},
  year         = {2024},
  url          = {https://arxiv.org/abs/2410.11523},
}

@misc{Verreault,
  author       = {Verreault, William},
  title        = {Almost sure upper bound for sums of random multiplicative
                  functions and critical chaos},
  howpublished = {arXiv:2608.21354},
  year         = {2026},
  url          = {https://arxiv.org/abs/2608.21354},
}

@misc{XuZheng2026TypicalSize,
  author       = {Xu, Max Wenqiang and Zheng, Junren},
  title        = {The typical size of {Hecke} eigenvalue sums is {$o(\sqrt{x})$}},
  howpublished = {arXiv:2609.14915},
  year         = {2026},
  url          = {https://arxiv.org/abs/2609.14915},
}

@article{ASVZZ,
  author  = {Aggarwal, Daksh and Subedi, Unique and Verreault, William
             and Zaman, Asif and Zheng, Chenghui},
  title   = {A conjectural asymptotic formula for multiplicative chaos in number theory},
  journal = {Res. Number Theory},
  volume  = {8},
  number  = {2},
  year    = {2022},
  pages   = {Paper No. 35},
}

@misc{OrsoniVerreault,
  author       = {Orsoni, Marcu-Antone and Verreault, William},
  title        = {Fourier decay of {Gaussian} multiplicative chaos and boundary geometry},
  year         = {2026},
  howpublished = {arXiv:2609.11761},
  url          = {https://arxiv.org/abs/2609.11761}
}

@article{SawinActa,
  author  = {Sawin, Will},
  title   = {Square-root cancellation for sums of factorization functions over squarefree progressions in {$\mathbb{F}_q[t]$}},
  journal = {Acta Math.},
  volume  = {233},
  number  = {2},
  year    = {2024},
  pages   = {285--418},
  doi     = {10.4310/ACTA.2024.v233.n2.a3}
}

@article {fusrp2,
    AUTHOR = {Hoban, Declan and Shah, Jibran Iqbal and Ismail,
              Nadya-Catherine and Verreault, William and Zaman, Asif},
     TITLE = {Central limit theorems for random multiplicative functions
              over function fields},
   JOURNAL = {Res. Math. Sci.},
  FJOURNAL = {Research in the Mathematical Sciences},
    VOLUME = {13},
      YEAR = {2026},
    NUMBER = {3},
     PAGES = {Paper No. 58, 38},}

@article {fusrp1,
    AUTHOR = {Hofmann, Maximilian C. E. and Hoganson, Annemily and Menon,
              Siddarth and Verreault, William and Zaman, Asif},
     TITLE = {Moments of random multiplicative functions over function
              fields},
   JOURNAL = {Math. Proc. Cambridge Philos. Soc.},
  FJOURNAL = {Mathematical Proceedings of the Cambridge Philosophical
              Society},
    VOLUME = {179},
      YEAR = {2025},
    NUMBER = {3},
     PAGES = {785--819},
}

\end{document}